\documentclass{article}
\usepackage{amssymb,amsmath}
\usepackage{graphics}
 
\def\qed{\hbox{${\vcenter{\vbox{                        
   \hrule height 0.4pt\hbox{\vrule width 0.4pt height 6pt
   \kern5pt\vrule width 0.4pt}\hrule height 0.4pt}}}$}}

\def\tr{\triangleright}

\newtheorem{theorem}{Theorem}

\newtheorem{lemma}[theorem]{Lemma}

\newtheorem{corollary}[theorem]{Corollary}

\newenvironment{proof}[1][Proof]{\smallskip\noindent{\bf #1.}\quad}%
{\qed\par\medskip}

\author{{\begin{tabular}{c} Gabriel Murillo \\
\small{\texttt{gmuri002@student.ucr.edu}}\end{tabular}}
\and
{\begin{tabular}{c} Sam Nelson \\
\small{\texttt{knots@esotericka.org}}\end{tabular}}
\and
\begin{tabular}{c}
\small{Department of Mathematics, University of California, Riverside} \\
\small{900 University Avenue, Riverside, CA, 92521}
\end{tabular}
}

\date{}
                           
\title{\Large \textbf{Alexander quandles of order 16}}

\begin{document}                                              
                                  
\maketitle

\abstract{
Isomorphism classes of Alexander quandles of order 16 are 
determined, and classes of connected quandles are identified. This paper 
extends the list of distinct connected finite Alexander quandles.}



\section{\large \textbf{Introduction}}

A \textit{quandle} is a set $Q$ with a non-associative binary
operation $\tr:Q\times Q\to Q$ satisfying
\newcounter{q}
\begin{list}{(\roman{q})}{\usecounter{q}}
\item for every $a\in Q$, we have $a\tr a= a$,
\item for every pair $a,b\in Q$ there is a unique $c\in Q$ such that 
$a=c\tr b$, and
\item for every $a,b,c\in Q$, we have $(a\tr b) \tr c = (a\tr c) \tr (b\tr c)$.
\end{list}

The three quandle axioms essentially form an algebraic distillation of the 
the three Reidemeister moves, which naturally makes quandles useful for
defining invariants of knots and links. In \cite{J}, the \textit{fundamental 
quandle} of a topological space is defined, and a Wirtinger-style 
presentation by generators and relations is given for the fundamental 
quandle of a knot or link complement. As with groups, distinguishing
quandles defined by generators and relations is a non-trivial problem
itself, but various techniques exist for using quandles to distinguish
knots, such as the 2-cocycle invariants defined in \cite{CJKLS}.

For the purpose of computing knot invariants using quandles, it is useful 
to compute isomorphism classes of finite quandles, particularly those of
finite connected quandles. Every group is a quandle with quandle operation 
given by conjugation, i.e. $a\tr b = b^{-1}ab$. Indeed, any union of 
conjugacy classes in a group forms a quandle.

One important class of quandles is the category of \textit{Alexander quandles}.
Let $M$ be any module over the ring $\Lambda=\mathbb{Z}[t^{\pm 1}]$ of Laurent 
polynomials in one variable. Then $M$ is a quandle with quandle operation
given by 
\[ a\tr b = ta + (1-t)b.\]

A quandle is \textit{connected} if it has a single orbit under $\tr$; for 
Alexander quandles, this is equivalent to $(1-t)M=M$. Connected quandles
are of particular importance in applications of quandle theory to knot theory
since all knot quandles are connected.

In \cite{N}, a method was given for determining all distinct isomorphism 
classes of Alexander quandles of a given finite order $n$, and the numbers
of distinct isomorphism classes were listed for values of $n$ up to 15.
In this paper, we compute all distinct isomorphism classes of Alexander 
quandles with 16 elements and identify which of these are connected.

\section{\large \textbf{Computations}}

An abelian group $M$ may be given the structure of a $\Lambda$-module, and
hence an Alexander quandle, by defining $tm=\phi(m)$ where 
$\phi\in\mathrm{Aut}_{\mathbb{Z}}(M)$ for each $m\in M$. Note that $\phi$ 
must be an automorphism in order to define multiplication by $t^{-1}$ as 
$\phi^{-1}$.

We will use the following theorem, proved in \cite{N}:

\begin{theorem}\label{tn}
Let $M$ and $M'$ be Alexander quandles with finite cardinality $|M|=|M'|$.
Then there is an isomorphism of Alexander quandles $f:M\to M'$ iff there is
an isomorphism of $\Lambda$-modules $h:(1-t)M\to (1-t)M'$.
\end{theorem}

That is, we can compare Alexander quandles by comparing their 
$\Lambda$-submodules $\mathrm{Im}(1-t)$. It is also useful to note the 
following lemma:

\begin{lemma} \label{conj}
Let $t_1, t_2\in \mathrm{Aut}_{\mathbb{Z}}(M)$ be $\mathbb{Z}$-automorphisms 
of $M$. Then the $\Lambda$-module structures $(M,t_1)$ and $(M,t_2)$ 
determined by $t_1$ and $t_2$ are isomorphic iff $t_1$ is conjugate to $t_2$.
\end{lemma}

\begin{proof}
Let $f:(M,t_1)\to (M,t_2)$ be an isomorphism of $\Lambda$-modules. Then
$f(t_1m)=t_2f(m)$ says that $t_1=f^{-1}t_2f\in \mathrm{Aut}_{\mathbb{Z}}(M)$.
Conversely, if there is an $f\in \mathrm{Aut}_{\mathbb{Z}}(M)$ with 
$t_1=f^{-1}t_2f$, then we have $f(t_1m)=t_2f(m)$ and $f$ is a 
$\mathbb{Z}$-automorphism which takes the action of $t_1$ to the action of 
$t_2$ on $M$; that is, $f$ is an isomorphism of $\Lambda$-modules.
\end{proof}

We may divide the problem of determining isomorphism classes of Alexander 
quandles with 16 elements into cases depending on which abelian group $M$ of
order 16 forms the underlying abelian group. The possibilities are
$M=\mathbb{Z}_{16}$ (the \textit{linear} Alexander quandles), 
$M=(\mathbb{Z}_2)^4$, $M=\mathbb{Z}_4\oplus \mathbb{Z}_4$, 
$M=\mathbb{Z}_8\oplus \mathbb{Z}_2$, and 
$M=\mathbb{Z}_4\oplus\mathbb{Z}_2\oplus\mathbb{Z}_2$.

\subsection{Linear Alexander quandles of order 16}

We are able to compute the linear Alexander quandles of order 16 by using
Corollary 2.2 in \cite{N}, which says:

\begin{corollary}
Let $a$ and $b$ be coprime to  $n$. Then the Alexander quandles 
$\Lambda_n/t-a$ and $\Lambda_n/t-b$ are isomorphic iff $N(n,a)=N(n,b)$ 
and $a\equiv b \mod N(n,a)$, where $N(n,a)=\frac{n}{\gcd(n,1-a)}$.
\end{corollary}

To do this we must first compute $N(16,a)=\frac{16}{\gcd(16,1-a)}$ for each 
$a\in \mathbb{Z}_{16}$ coprime to $16$. This yields four different values 
for $N$; from the corollary we can tell there are three pairs of linear 
Alexander quandles which are isomorphic to one another, namely 
$\Lambda_{16}/t-3 \cong \Lambda_{16}/t-11$, $\Lambda_{16}/t-7 \cong 
\Lambda_{16}/t-15$, and $\Lambda_{16}/t-5 \cong \Lambda_{16}/t-13$. The 
quandle $\Lambda_{16}/t-9$ and the trivial quandle 
$\Lambda_{16}/t-1\cong T_{16}$ form their own isomorphism classes.

For the purpose of comparing the linear Alexander quandles of order 16 with
other Alexander quandles of order 16, we must still compute the submodules
$\mathrm{Im}(1-t)$ for each linear quandle. The result are summarized in 
figure \ref{lin}.

\begin{figure} [!ht] \label{lin}
\begin{center}
\begin{tabular}{|c|c|} \hline
$M$  & $(1-t)M$ \\ \hline
$\Lambda_{16}/t-1$  & $0$ \\
$\Lambda_{16}/t-3$  & $\Lambda_8/t-3$ \\
$\Lambda_{16}/t-11$ & $\Lambda_8/t-3$ \\ 
$\Lambda_{16}/t-7$  & $\Lambda_8/t-7$ \\
$\Lambda_{16}/t-15$ & $\Lambda_8/t-7$ \\ 
$\Lambda_{16}/t-5$  & $\Lambda_4/t+3$ \\ 
$\Lambda_{16}/t-13$ & $\Lambda_4/t+3$ \\  
$\Lambda_{16}/t-9$  & $\Lambda_2/t+1$ \\ \hline
\end{tabular}
\end{center}
\caption{$\mathrm{Im}(1-t)$ for linear Alexander quandles of order 16.}
\end{figure}

\subsection{Alexander quandles structures on $(\mathbb{Z}_2)^4$}

To compute this case, we realize that $(\mathbb{Z}_2)^4$ is not just a
$\mathbb{Z}$-module but also a $\mathbb{Z}_2$-module, and since
$\mathbb{Z}_2$ is a principal ideal domain, we are able to use
the classification theorem for modules over a PID. We start by listing
all $\Lambda$-modules whose underlying abelian group is $(\mathbb{Z}_2)^4$. 
Specifically, these have the form 
\[\bigoplus_{i=1}^n \Lambda_2/h_i, \qquad \mathrm{where} \ \ 
h_1|h_2|\dots |h_n, \quad\mathrm{and} \quad \sum_{i=1}^n \mathrm{deg}(h_i)=4.\]
With this list and theorem \ref{tn}, we are able to compare the Alexander 
quandles by computing $(1-t)/M$ for each module $M=\Lambda_2/h$ where 
$h\in \Lambda_2$ is a polynomial over $\mathbb{Z}_2$ with lead coefficient 1
and constant term 1. That is, we must multiply every
element in each module $M=\Lambda_2/h$ by $(1-t)$ and reduce modulo $h$.
Then we must identify the resulting submodule. Note that when 
$h(a,b,c)=1+at+bt^2+ct^3+t^4$ with $a,b,c\in \mathbb{Z}_2$  has an even 
number of terms, the submodules $(1-t)M/h(a,b,c)$ are equal as sets (though 
distinct as $\Lambda$-modules), and the same is true for all $h(a,b,c)$ with 
an odd number of terms.

\begin{figure}[!ht]
\begin{center}
\begin{tabular}{|c|c|c|} \hline
Connected & $M$ & $(1-t)M$ \\ \hline
* & $\Lambda_2/t^4+t+1$ & $\Lambda_2/t^4+t+1$ \\
* & $\Lambda_2/t^4+t^2+1$ & $\Lambda_2/t^4+t^2+1$ \\
* & $\Lambda_2/t^4+t^3+1$ & $\Lambda_2/t^4+t^3+1$ \\
* & $\Lambda_2/t^4+t^3+t^2+t+1$ & $\Lambda_2/t^4+t^3+t^2+t+1$ \\ \hline
& $\Lambda_2/t^4+1$ & $\Lambda_2/t^3+t^2+t+1$  \\
& $\Lambda_2/t^4+t^2+t+1$ & $\Lambda_2/t^3+t^2+1$  \\
& $\Lambda_2/t^4+t^3+t^2+1$ & $\Lambda_2/t^3+t+1$  \\
&  $\Lambda_2/t^4+t^3+t+1$ & $\Lambda_2/t^3+1$   \\  \hline
&  $(\Lambda_2/t+1)^4$ & $0$ \\
&  $(\Lambda_2/t+1)^2\oplus\Lambda_2/t^2+1$ & $\Lambda_2/t+1$ \\
&  $(\Lambda_2/t^2+1)^2$ & $(\Lambda_2/t+1)^2$ \\
* &  $(\Lambda_2/t^2+t+1)^2$ & $(\Lambda_2/t^2+t+1)^2$ \\
&  $\Lambda_2/t+1\oplus\Lambda_2/t^3+1$ & $\Lambda_2/t^2+t+1$ \\
&  $\Lambda_2/t+1\oplus\Lambda_2/t^3+t^2+t+1$ & $\Lambda_2/t^2+1$  \\ \hline
\end{tabular}
\end{center}
\caption{$\mathrm{Im}(1-t)$ for Alexander quandles with abelian group 
$(\mathbb{Z}_2)^4$ }
\end{figure}

\subsection{Alexander quandles defined by $\mathbb{Z}$-automorphisms}

In light of lemma \ref{conj}, it is sufficient to consider only a single
representative from each conjugacy class. Using a maple program, we first
compute the $\mathbb{Z}$-automorphism group of $M$. We represent an element
$\phi\in\mathrm{Aut}_{\mathbb{Z}}(M)$ by listing an image for each element of 
a generating set. The program checks each such set of images to determine 
whether the linear map thus defined is an automorphism. The program then 
compares the automorphisms pairwise for conjugacy and deletes redundant 
automorphisms, yielding a single representative for each conjugacy class.
We then compute the $\Lambda$-submodule $\mathrm{Im}(\mathrm{Id}-\phi)$
for each representative automorphism $\phi$. 

We applied this procedure for $\mathbb{Z}_8\oplus\mathbb{Z}_2$ and
$\mathbb{Z}_4\oplus\mathbb{Z}_2\oplus\mathbb{Z}_2$. The results are collected 
in table \ref{res}. Two of the classes of linear quandles are isomorphic to
quandles listed in the table, namely $\Lambda_{16}/t-9\cong 
(\Lambda_2/t+1)^2\oplus \Lambda_2/t^2+1$ and $\Lambda_{16}/t-5\cong 
(\mathbb{Z}_4\oplus\mathbb{Z}_4, \phi(1,0)=(0,1),\phi(0,1)=(3,2))$.

Together with the results from the previous two sections, we have our main
result, namely:
\begin{theorem}
There are a total of 23 distinct isomorphism classes of Alexander quandles 
with 16 elements. Of these, eight are connected, including five quandles with 
underlying abelian group $(\mathbb{Z}_2)^4$ and three quandles with 
underlying abelian group $\mathbb{Z}_4\oplus\mathbb{Z}_4$.
\end{theorem}

\begin{proof}
In light of theorem \ref{tn}, this is simply a matter of counting distinct
submodules $\mathrm{Im}(1-t)$. In all, a total of 23 
distinct submodules appear; of these, eight are connected, namely
$\Lambda_2/t^4+t+1$, $\Lambda_2/t^4+t^2+1$, $\Lambda_2/t^4+t^3+1$, 
$\Lambda_2/t^4+t^3+t^2+t+1$, $(\Lambda_2/t^2+t+1)^2$ and 
$\mathbb{Z}_4\oplus\mathbb{Z}_4$ with Alexander quandle structure given by 
$\phi(1,0)=(1,0)$ and $\phi(0,1)=(1,1), (3,1)$ and $(3,3)$. 
\end{proof}

\begin{figure}[!ht] \label{res}

\begin{center}
\begin{tabular}{cc}
\begin{tabular}{|r|cc|c|} 
\multicolumn{4}{c}{Automorphisms of $\mathbb{Z}_4\oplus\mathbb{Z}_4$} \\
\multicolumn{4}{c}{}  \\ \hline
& $\phi((1,0))$ & $\phi((0,1))$ & $\mathrm{Im}(\mathrm{Id}-\phi)$  \\ \hline
& $(1,0)$ & $(0,1)$ & $0$ \\ \hline
& $(1,0)$ & $(0,3)$ & $\Lambda_2/t+1$ \\
& $(1,0)$ & $(2,1)$ & $\Lambda_2/t+1$ \\
& $(1,2)$ & $(2,1)$ & $\Lambda_2/t+1$ \\ \hline
& $(1,2)$ & $(2,3)$ & $(\Lambda_2/t+1)^2$ \\
& $(3,0)$ & $(0,3)$ & $(\Lambda_2/t+1)^2$ \\
& $(0,1)$ & $(3,2)$ & $\Lambda_4/t+3$ \\
& $(0,1)$ & $(1,0)$ & $\Lambda_4/t+1$ \\ \hline
* & $(0,1)$ & $(1,1)$ & $(\mathbb{Z}_4\oplus\mathbb{Z}_4, \phi)$ \\
* & $(0,1)$ & $(3,1)$ & $(\mathbb{Z}_4\oplus\mathbb{Z}_4, \phi)$ \\
* & $(0,1)$ & $(3,3)$ & $(\mathbb{Z}_4\oplus\mathbb{Z}_4, \phi)$ \\ \hline
& $(0,1)$ & $(3,0)$ & $(\mathbb{Z}_4\oplus\mathbb{Z}_2, \phi')$  \\ 
& $(0,1)$ & $(1,2)$ & $(\mathbb{Z}_4\oplus\mathbb{Z}_2, \phi'')$ \\  \hline

\multicolumn{4}{c}{} \\
\multicolumn{4}{c}{
$\phi'(1,0)=(3,1), \phi'(0,1)=(2,1)$}  \\
\multicolumn{4}{c}{ 
$\phi''(1,0)=(1,1), \phi''(0,1)=(2,1)$} \\
\end{tabular} 
&
\begin{tabular}{|cc|c|} 
\multicolumn{3}{c}{Automorphisms of $\mathbb{Z}_8\oplus\mathbb{Z}_2$} \\
\multicolumn{3}{c}{}  \\ \hline
$\phi((1,0))$ & $\phi((0,1))$ & $\mathrm{Im}(\mathrm{Id}-\phi)$  \\ \hline
$(1,0)$ & $(0,1)$ & $0$ \\ \hline
$(1,1)$ & $(0,1)$ & $\Lambda_2/t+1$ \\
$(1,0)$ & $(4,1)$ & $\Lambda_2/t+1$ \\
$(5,0)$ & $(0,1)$ & $\Lambda_2/t+1$ \\ \hline
$(1,1)$ & $(4,1)$ & $\Lambda_2/t^2+1$ \\ \hline
$(3,0)$ & $(0,1)$ & $\Lambda_4/t+1$ \\
$(3,0)$ & $(4,1)$ & $\Lambda_4/t+1$ \\
$(3,1)$ & $(0,1)$ & $\Lambda_4/t+1$ \\
$(7,0)$ & $(0,1)$ & $\Lambda_4/t+1$ \\ \hline
$(3,1)$ & $(4,1)$ & $\Lambda_4/t+3$ \\ \hline
\end{tabular}  \\
\end{tabular}

\medskip

\begin{tabular}{|ccc|c|}
\multicolumn{4}{c}{Automorphisms of 
$\mathbb{Z}_4\oplus\mathbb{Z}_2\oplus\mathbb{Z}_2$} \\
\multicolumn{4}{c}{}  \\ \hline
$\phi((1,0,0))$ & $\phi((0,1,0))$ & $\phi((0,0,1))$ & 
$\mathrm{Im}(\mathrm{Id}-\phi)$  \\ \hline
$(1,0,0)$ & $(0,1,0)$ & $(0,0,1)$ & $0$ \\ \hline
$(1,0,0)$ & $(0,0,1)$ & $(0,1,0)$ & $\Lambda_2/t+1$ \\
$(1,0,0)$ & $(0,1,0)$ & $(2,0,1)$ & $\Lambda_2/t+1$ \\
$(1,0,1)$ & $(0,1,0)$ & $(0,0,1)$ & $\Lambda_2/t+1$ \\
$(3,0,0)$ & $(0,1,0)$ & $(0,0,1)$ & $\Lambda_2/t+1$ \\ \hline
$(1,0,1)$ & $(2,1,0)$ & $(0,0,1)$ & $(\Lambda_2/t+1)^2$ \\
$(1,0,1)$ & $(2,1,1)$ & $(0,0,1)$ & $(\Lambda_2/t+1)^2$ \\ \hline
$(1,0,0)$ & $(0,0,1)$ & $(2,1,0)$ & $\Lambda_2/t^2+1$ \\
$(1,0,1)$ & $(0,0,1)$ & $(0,1,0)$ & $\Lambda_2/t^2+1$ \\
$(1,0,1)$ & $(0,1,0)$ & $(2,0,1)$ & $\Lambda_2/t^2+1$ \\ 
$(1,0,0)$ & $(0,0,1)$ & $(0,1,1)$ & $\Lambda_2/t^2+t+1$ \\ \hline 
$(1,0,1)$ & $(0,0,1)$ & $(2,1,1)$ & $\Lambda_2/t^3+1$ \\
$(1,0,1)$ & $(0,0,1)$ & $(2,1,0)$ & $\Lambda_2/t^3+t^2+t+1$ \\ \hline
\end{tabular}
\end{center}

\noindent * Connected quandles.

\caption{Results of computation of $\mathrm{Im}(\mathrm{Id}-\phi)$ for 
Alexander quandles given by automorphisms.}
\end{figure}

\end{document}